\documentstyle[amssymb,12pt]{amsart}

\newtheorem{theorem}{Theorem}
\newtheorem{lemma}[theorem]{Lemma}

\newcommand{\forces}{\Vdash}
\newcommand{\res}{\upharpoonright}
\newcommand{\ch}{$C\!H$ and $2^{\omega_1}\!>\!\omega_2\;$}

\begin{document}

\baselineskip=18pt

  \begin{center}
     {\large Essential Kurepa Trees Versus\\ 
     Essential Jech--Kunen Trees}
  \end{center}

  \begin{center}
 Renling Jin\footnote{The first author 
would like to thank Mathematics Department of 
Rutgers University for its hospitality during his one week visit there 
in October 1992, when the main part of the paper was developed.}
\& Saharon Shelah\footnote{Publ. 
nu. 498. The second
author was partially supported by the United States--Israel Binational
Science Foundation.}
  \end{center}

  \bigskip

  \begin{quote}

    \centerline{Abstract}

    \small

By an $\omega_{1}$--tree we mean a tree of size $\omega_{1}$ and height
$\omega_{1}$. An $\omega_1$--tree is called a Kurepa tree if all its levels
are countable and it has more than $\omega_1$ branches.
An $\omega_{1}$--tree is called a Jech--Kunen tree if it has
$\kappa$ branches for some $\kappa$ strictly between $\omega_{1}$
and $2^{\omega_{1}}$. A Kurepa tree is called an essential Kurepa tree
if it contains no Jech--Kunen subtrees. 
A Jech--Kunen tree is called an essential Jech--Kunen tree 
if it contains no Kurepa subtrees. In this paper
we prove that (1) it is consistent with $C\!H$ and $2^{\omega_1}>\omega_2$ that
there exist essential Kurepa trees and there are no essential Jech--Kunen 
trees, (2) it is consistent with $C\!H$ and $2^{\omega_1}>\omega_2$ plus
the existence of a Kurepa tree with $2^{\omega_1}$ branches that 
there exist essential Jech--Kunen trees and there are no essential Kurepa 
trees. In the second result we require the existence of a Kurepa tree with 
$2^{\omega_1}$ branches in order to avoid triviality.

\end{quote}

\section{Introduction}

Our trees are always growing downward. We use $T_{\alpha}$ 
for the $\alpha^{th}$ level of $T$ and use $T\!\res\!\alpha$ for
$\bigcup_{\beta<\alpha}T_{\beta}$. For every $t\in T$ let $ht(t)=\alpha$
iff $t\in T_{\alpha}$. Let $ht(T)$, the height of $T$, be the least
ordinal $\alpha$ such that $T_{\alpha}=\emptyset$. By a branch of $T$
we mean a totally ordered subset of $T$ which intersects every nonempty
level of $T$. For any tree $T$ let $m(T)$ be the set of all maximal
nodes of $T$, {\em i.e.} $m(T)=\{t\in T:(\forall s\in T)(s\leqslant t\rightarrow
s=t)\}$. All trees considered in this paper have cardinalities 
less than or equal to $\omega_1$ so that, without loss of generality,
we can assume all those trees are subtrees of $(\omega_1^{<\omega_1},
\supseteq)$, where $\omega_1^{<\omega_1}$ is the set of all functions
from some countable ordinals to $\omega_1$. Hence every tree here
has a unique root $\emptyset$ and if 
$\{t_n:n\in\omega\}\subseteq T$
is a decreasing sequence of $T$, then $t=\bigcup_{n\in\omega}t_n$ is the only
possible greatest lower bound of $\{t_n:n\in\omega\}$. We are also free
to use either $\leqslant_T$ or $\supseteq$ for the order of a tree $T$,
{\em i.e.} $s\leqslant_T t$ if and only if $s\supseteq t$.

By an $\omega_1$--tree we mean a tree of height $\omega_1$ and size
$\omega_1$. Notice that our definition of $\omega_1$--tree is slightly
different from the usual definition by not 
requiring every level to be countable.
An $\omega_1$--tree $T$ is called a {\em Kurepa tree} if every level of $T$ is
countable and $T$ has more than $\omega_1$ branches. An $\omega_1$--tree
$T$ is called a {\em Jech--Kunen tree} if $T$ has $\kappa$ branches for
some $\kappa$ strictly between $\omega_1$ and $2^{\omega_1}$. 
We call a Kurepa tree {\em thick} if it has $2^{\omega_1}$ branches.
Obviously, a Kurepa non-Jech--Kunen tree must be thick, and a Jech--Kunen
tree with every level countable is a Kurepa tree.

While Kurepa trees are better studied, Jech--Kunen trees are relatively less
popular. It is K. Kunen [K1][Ju], who brought Jech--Kunen trees to people's
attention by proving that: under $C\!H$ and $2^{\omega_1}\!>\omega_2$, 
the existence of a compact Hausdorff space with weight $\omega_1$ and
size strictly between $\omega_1$ and $2^{\omega_1}$ is equivalent
to the existence of a Jech--Kunen tree. It is also easy to observe that:
under $C\!H$ and $2^{\omega_1}\!>\omega_2$, the existence of a (Dedekind)
complete dense linear order with density $\omega_1$ and size
strictly between $\omega_1$ and $2^{\omega_1}$ is also equivalent to
the existence of a Jech--Kunen tree. Above results are interesting
because those compact Hausdorff spaces and complete dense linear orders
cannot exist if we replace $\omega_1$ by $\omega$, while the existence
of a Jech--Kunen tree is undecidable. In this paper we would like to
consider Jech--Kunen trees only under \ch.

The consistency of a Jech--Kunen tree was given in [Je1], in which T. Jech
constructed a generic Kurepa tree with less than $2^{\omega_1}$ branches
in a model of $C\!H$ and $2^{\omega_1}\!>\omega_2$. By assuming the consistency
of an inaccessible cardinal, K. Kunen proved the consistency of 
non--existence of Jech--Kunen trees with 
\ch (see [Ju, Theorem 4.8]). In Kunen's model
there are also no Kurepa trees. Kunen proved (see [Ju, Theorem 4.10])
also that the assumption of an inaccessible cardinal above is necessary.
The differences between Kurepa trees and 
Jech--Kunen trees in terms of the existence have been studied in 
[Ji1] [Ji2] [Ji3] [SJ1] [SJ2]. It was proved that
the consistency of an inaccessible cardinal implies (1) it is consistent
with \ch that there exist Kurepa trees but there are no Jech--Kunen trees 
[SJ1], (2) it is consistent with \ch that there exist Jech--Kunen trees but 
there are no Kurepa trees [SJ2].

What could we say without the presence of large cardinals?
In stead of killing all Kurepa trees, which needs an inaccessible
cardinal, while keeping some Jech--Kunen trees alive, or killing all
Jech--Kunen trees, which needs again an inaccessible cardinal, while
keeping some Kurepa trees alive,
we can kill all Kurepa subtrees of a Jech--Kunen tree or kill
all Jech--Kunen subtrees of a Kurepa tree without using large cardinals.
Let's call a Kurepa tree $T$ essential if $T$ has no Jech--Kunen subtrees,
and call a Jech--Kunen tree $T$ essential if $T$ has no Kurepa subtrees.
In [Ji1], the first author proved that it is consistent with \ch,
together with {\em Generalized Martin's Axiom} and the existence of a thick
Kurepa tree, that no essential Kurepa trees
and no essential Jech--Kunen trees.
We required the presence of thick Kurepa trees in the model in order to
avoid triviality. In [Ji3], the first author proved that it is
consistent with \ch that there exist both essential Kurepa trees and 
essential Jech--Kunen trees. A weak version of this result was proved 
in [Ji1] with help of an inaccessible cardinal. This paper could be
considered as a continuation of the research
done in [Ji1] [Ji2] [Ji3] [SJ1] [SJ2].

In \S 1, we prove that it is consistent with \ch that there exist
essential Kurepa trees but there are no essential Jech--Kunen trees.
In \S 2, we prove that it is consistent with \ch plus the existence
of a thick Kurepa tree that there exist
essential Jech--Kunen trees but there are no essential Kurepa trees.
In \S 3, we simplify the proofs of two old results by using the forcing notion
for producing a generic essential Jech--Kunen tree defined in \S 2.

We write $\dot{a}$ in the ground model for a name of an element $a$ in 
the forcing extension. If $a$ is in the ground model, we usually write
$a$ itself as a canonical name of $a$. The rest of the notation will
be consistent with [K2] or [Je2].

\section{Yes Essential Kurepa Trees, No Essential Jech--Kunen Trees.}

In this section we are going to construct a model of \ch in which there
exist essential Kurepa trees and there are no essential Jech--Kunen trees.
Our strategy to do this can be described as follows: 
first, we take a model of \ch
plus $G\!M\!A$ (Generalized Martin's Axiom) as our ground model, so that
in the ground model there are no essential Jech--Kunen trees, then, we
add a generic Kurepa tree which has no Jech--Kunen subtrees. 
The hard part is to prove that the forcing adds no essential 
Jech--Kunen trees.

Let $\Bbb P$ is a poset. A subset $S$ of $\Bbb P$ is called {\em linked}
if any two elements in $S$ is compatible in $\Bbb P$. A poset $\Bbb P$ is
called $\omega_1$--{\em linked} if $\Bbb P$ is the union of 
$\omega_1$ linked subsets of $\Bbb P$. A subset $S$ of $\Bbb P$ is called 
{\em centered} if every finite subset of $S$ has a lower bound in $\Bbb P$.
A poset $\Bbb P$ is called {\em countably compact} if every countable
centered subset of $\Bbb P$ has a lower bound in $\Bbb P$. Now $G\!M\!A$
is the following statement:

\begin{quote}
Suppose $\Bbb P$ is an $\omega_1$--linked and countably compact poset.
For any $\kappa<2^{\omega_1}$, if ${\cal D}=\{D_{\alpha}:\alpha<\kappa\}$
is a collection of $\kappa$ dense subsets of $\Bbb P$, then there
exists a filter $G$ of $\Bbb P$ such that $G\cap D_{\alpha}\not=\emptyset$
for all $\alpha<\kappa$.
\end{quote}

We choose the form of $G\!M\!A$ from [B], where a model of \ch plus
$G\!M\!A$ can be found.

Let $I$ be any index set. We write ${\Bbb K}_{I}$ for a poset such that
$p$ is a condition in ${\Bbb K}_{I}$ iff 
$p=(A_p,l_p)$ where $A_p$ is a countable
subtree of $(\omega_1^{<\omega_1},\supseteq)$ of height $\alpha_p+1$ and
$l_p$ is a function from a countable subset of $I$ into
$(A_p)_{\alpha_p}$, the top level of $A_p$. For any $p,q\in
{\Bbb K}_I$, define $p\leqslant q$ iff 

(1) $A_p\!\res\!\alpha_q+1=A_q$, 

(2) $dom(l_p)\supseteq dom(l_q)$ 

(3) $(\forall\xi\in dom(l_q))(l_q (\xi) \subseteq l_p (\xi))$.

It is easy to see that ${\Bbb K}_I$ is countably closed (or $\omega_1$--closed).
If $C\!H$ holds, then ${\Bbb K}_I$ is $\omega_1$--linked. Let $M$ be a model
of $C\!H$ and ${\Bbb K}_I\in M$. Suppose that $G$ is a ${\Bbb K}_I$--generic
filter over $M$ and let $T_G=\bigcup_{p\in G}A_p$. Then in $M[G]$, the tree
$T_G$ is an $\omega_1$--tree with every level countable and $T_G$ has
exactly $|I|$ branches. Furthermore, if for every $i\in I$ let
\[B(i)=\bigcup\{l_p(i):p\in G\mbox{ and }i\in dom(l_p)\},\]
then $B(i)\not=B(i')$ for any $i, i'\in I$ and $i\not=i'$, and  
$\{B(i):i\in I\}$ is the set of all branches of $T_G$ in $M[G]$.
Hence if $|I|>\omega_1$,
then $T_G$ will be a Kurepa tree with $|I|$ branches in $M[G]$.
${\Bbb K}_I$ is the poset used in [Je1] for creating a generic Kurepa tree.
All those facts above can also be found in [Je1] or [T].

For convenience we sometimes view ${\Bbb K}_I$ as an iterated forcing
notion \[{\Bbb K}_{I'}*Fn(I\smallsetminus I',T_{\dot{G}_{I'}},\omega_1),\]
for any $I'\subseteq I$, where $G_{I'}$ is a ${\Bbb K}_{I'}$--generic
filter over the ground model and $Fn(I\smallsetminus 
I',T_{G_{I'}},\omega_1)$, in $M[G_{I'}]$, is the set of all functions 
from some countable subset of $I\smallsetminus I'$ to 
$T_{G_{I'}}$ with the order defined by letting $p\leqslant q$ iff
$dom(q)\subseteq dom(p)$ and for any $i\in dom(q)$, $p(i)\leqslant q(i)$. 
The poset $Fn(J,T_{G},\omega_1)$ is in fact the countable support product
of $|J|$--copies of $T_{G}$. We call two posets $\Bbb P$ and $\Bbb Q$ are
forcing equivalent if there is a poset $\Bbb R$ such that $\Bbb R$ can be
densely embedded into both $\Bbb P$ and $\Bbb Q$. The posets ${\Bbb K}_I$ and 
${\Bbb K}_{I'}*Fn(I\smallsetminus I',T_{\dot{G}_{I'}},\omega_1)$
are forcing equivalent because the map \[F:{\Bbb K}_I\mapsto
{\Bbb K}_{I'}*Fn(I\smallsetminus I',T_{\dot{G}_{I'}},\omega_1)\]
such that for every $p\in {\Bbb K}_{I}$, \[F(p)=((A_p,l_p\!\res\! I'),
l_p\!\res\! I\smallsetminus I')\]  
is a dense embedding. 

\begin{lemma}[K. Kunen]

Let $M$ be a model of $C\!H$. Suppose that $\lambda>\omega_2$ is a cardinal
in $M$ and ${\Bbb K}_{\lambda}\in M$. Suppose $G_{\lambda}$ is a 
${\Bbb K}_{\lambda}$--generic filter over $M$ and $T_{G_{\lambda}}=
\bigcup_{p\in G_{\lambda}}A_p$. Then in $M[G_{\lambda}]$ 
the tree $T_{G_{\lambda}}$ is a Kurepa tree with $\lambda$ 
branches and $T_{G_{\lambda}}$ has no subtrees with
$\kappa$ branches for some $\kappa$ strictly between $\omega_1$ and 
$\lambda$.

\end{lemma}

\noindent {\bf Proof:}\quad
Assume that $T$ is a subtree of $T_{G_{\lambda}}$ with
more than $\omega_1$ branches in $M[G_{\lambda}]$. We want to show that
$T$ has $\lambda$ branches in $M[G_{\lambda}]$. 
Since $|T|=\omega_1$, then there exists a subset $I\subseteq \lambda$
in $M$ with cardinality $\leqslant\omega_1$ such that
$T\in M[G_{I}]$, where \[G_I=\{p\in G:dom(l_p)\subseteq I\}.\] 
Notice that $T_{G_{\lambda}}=T_{G_I}$ (in fact $T_G=T_{G_{\emptyset}}$).
Since in $M[G_{I}]$ the tree $T_{G_I}$ has only $|I|$ branches, then 
the tree $T$ can have at most $\omega_1$ branches in $M[G_I]$. 
Let $B$ be a branch of $T$ in $M[G_{\lambda}]$ which is not in
$M[G_I]$. Since $|B|=\omega_1$, there exists a subset
$J$ of $\lambda\smallsetminus I$ with cardinality $\leqslant\omega_1$ such
that $B\in M[G_I][H_J]$, where $H_J$ is a $Fn(J,T_{G_I},\omega_1)$--generic 
filter over $M[G_I]$. Now $\lambda\smallsetminus I$
can be partitioned into $\lambda$--many subsets of cardinality $\omega_1$
and for every subset $J'\subseteq\lambda\smallsetminus (I\cup J)$ 
of cardinality $\omega_1$
the poset ${\Bbb P}_J=Fn(J,T_{G_I},\omega_1)$ is isomorphic to
the poset ${\Bbb P}_{J'}=Fn(J',T_{G_I},\omega_1)$ through
an obvious isomorphism $\pi$ induced by a bijection between $J$ and $J'$. 
Let $\dot{B}$ be a ${\Bbb P}_J$--name for $B$. Then $\pi_*(\dot{B})$
is a ${\Bbb P}_{J'}$--name for a new branch of $T$. 
Forcing with ${\Bbb P}_J\times {\Bbb P}_{J'}$ will create two different
branches $\dot{B}_{H_J}$ and $(\pi_*(\dot{B}))_{H_{J'}}$. Hence
forcing with $Fn(\lambda\smallsetminus I,T_{G_I},\omega_1)$ 
will produce at least $\lambda$ new branches of $T$.
\quad $\square$ 

\bigskip

Next lemma is a simple fact which will be used later.

\begin{lemma}

Suppose ${\Bbb P}$ is an $\omega_1$--closed poset of size $\omega_1$
(hence $C\!H$ must hold). Then the tree $(\omega_1^{<\omega_1},\supseteq)$
can be densely embedded into $\Bbb P$.

\end{lemma}

\noindent {\bf Proof:}\quad Folklore. \quad $\square$

\begin{lemma}

Let $M$ be a model of \ch plus $G\!M\!A$ and let ${\Bbb P}=
(\omega_1^{<\omega_1},\supseteq)\in M$. Suppose $G$ is a $\Bbb P$--generic
filter over $M$. Then in $M[G]$ every Jech--Kunen tree has a Kurepa subtree.

\end{lemma}

\noindent {\bf Proof:}\quad
Let $T$ be a Jech--Kunen tree in $M[G]$ with $\delta$ branches for 
$\omega_1<\delta<\lambda=2^{\omega_1}$. Without loss of generality 
we can assume that there is a regular cardinal $\kappa$ such that 
$\omega_1<\kappa\leqslant\delta$ and for every $t\in T$ there are at least 
$\kappa$ branches of $T$ passing through $t$ in $M[G]$. Again in $M[G]$
let $f:\kappa\mapsto {\cal B}(T)$ be a one to one function such that
for every $t\in T$ and for every $\alpha<\kappa$ there exists an $\beta\in
\kappa\smallsetminus\alpha$ such that $t\in f(\beta)$.
Without loss of generality let us assume that
\[1_{\Bbb P}\forces (\dot{T}\mbox{ is a Jech--Kunen tree and }
\dot{f}:\kappa\mapsto {\cal B}(\dot{T})\]\[\mbox{ is a one to one function
such that }(\forall t\in\dot{T})(\forall\alpha\in\kappa)(\exists\beta\in
\kappa\smallsetminus\alpha)(t\in\dot{f}(\beta))).\]
We want now to construct a poset $\Bbb R$ in $M$ such that a filter
$H$ of $\Bbb R$ obtained by applying $G\!M\!A$ in $M$ will
give us a $\Bbb P$--name for a Kurepa subtree of $T$ in $M[G]$.

Let $r$ be a condition of $\Bbb R$ iff $r=(I_r,{\Bbb P}_r,{\cal A}_r,
{\cal S}_r)$ where $I_r$ is a countable subtree of $(\omega_1^{<\omega_1},
\supseteq)$, ${\Bbb P}_r=\langle p^r_t:t\in I_r\rangle$, ${\cal A}_r=
\langle A^r_t:t\in I_r\rangle$ and ${\cal S}_r=\langle S^r_t:t\in I_r\rangle$
such that

(1) ${\Bbb P}_r\subseteq {\Bbb P}$ and for every $t\in I_r$ the element
$A^r_t$ is a nonempty countable subtree of 
$(\omega_1^{<\omega_1},\supseteq)$ of
height $\alpha^r_t+1$ (we will use some $A^r_t$'s to generate a 
Kurepa subtree of $T$) and $S^r_t$ is a nonempty
countable subset of $\kappa$,

(2) $(\forall s,t\in I_r)(s\subseteq t\leftrightarrow p^r_t\leqslant p^r_s)$,
(This implies that $s$ and $t$ are incompatible iff $p^r_s$ and $p^r_t$
are incompatible for $s,t\in I_r$ because $\Bbb P$ is a tree,)

(3) $(\forall s,t\in I_r)(s\subseteq t\rightarrow 
A^r_t\!\res\! ht(A^r_s)=A^r_s)$, 

(4) $(\forall s,t\in I_r)(s\subseteq t\rightarrow S^r_s\subseteq S^r_t)$,

(5) $(\forall t\in I_r)(p^r_t\forces A^r_t\subseteq\dot{T})$,

(6) $(\forall t\in I_r)(\forall\alpha\in S^r_t)(\exists a\in 
(A^r_t)_{\alpha^r_t})(p^r_t\forces a\in\dot{f}(\alpha))$.

The order of $\Bbb R$: for any $r,r'\in {\Bbb R}$, let
$r\leqslant r'$ iff $I_{r'}\subseteq I_r$ and for every $t\in I_{r'}$
\[p^{r'}_t=p^r_t,\; A^{r'}_t=A^r_t\; \mbox{ and }S^{r'}_t\subseteq S^r_t.\]

\bigskip

{\bf Claim 3.1}\quad The poset $\Bbb R$ is $\omega_1$--linked.

Proof of Claim 3.1:\quad Let $r,r'\in {\Bbb R}$ such that
$I_r=I_{r'}$, ${\Bbb P}_r={\Bbb P}_{r'}$ and ${\cal A}_r={\cal A}_{r'}$.
Then the condition $r''\in {\Bbb R}$ such that \[I_{r''}=I_r,\;
{\Bbb P}_{r''}={\Bbb P}_r,\; {\cal A}_{r''}={\cal A}_r\mbox{ and }
{\cal S}_{r''}=\langle S^r_t\cup S^{r'}_t:t\in I_{r''}\rangle\] is a common
lower bound of both $r$ and $r'$. Since there are only $\omega_1$ different
$\langle I_r,{\Bbb P}_r,{\cal A}_r\rangle$'s and for each fixed $\langle
I_{r_0},{\Bbb P}_{r_0},{\cal A}_{r_0}\rangle$ the set \[\{r\in {\Bbb R}:
\langle I_r,{\Bbb P}_r,{\cal A}_r\rangle=\langle I_{r_0},{\Bbb P}_{r_0},
{\cal A}_{r_0}\rangle\}\] is linked, then $\Bbb R$ is the union of
$\omega_1$ linked subsets of $\Bbb R$.\quad $\square$ (Claim 3.1)

\bigskip

{\bf Claim 3.2}\quad The poset ${\Bbb R}$ is countably compact.

Proof of Claim 3.2:\quad Suppose that ${\Bbb R}'$ is a countable
centered subset of $\Bbb R$. Notice that for any finite
${\Bbb R}'_0\subseteq {\Bbb R}'$ and for any $t\in\bigcap\{I_r:
r\in {\Bbb R}'_0\}$ all $p^r_t$'s are same and all $A^r_t$ are same 
for $r\in {\Bbb R}'_0$ because ${\Bbb R}'_0$ has a common lower bound
in $\Bbb R$. We now want to construct a condition $\bar{r}\in {\Bbb R}$
such that $\bar{r}$ is a common lower bound of ${\Bbb R}'$. Let

(1) $I_{\bar{r}}=\bigcup_{r\in {\Bbb R}'}I_r$,

(2) ${\Bbb P}_{\bar{r}}=\langle p^{\bar{r}}_t:t\in I_{\bar{r}}\rangle$
where $p^{\bar{r}}_t=p^r_t$ for some $r\in {\Bbb R}'$ such that
$t\in I_r$,

(3) ${\cal A}_{\bar{r}}=\langle A^{\bar{r}}_t:t\in I_{\bar{r}}\rangle$
where $A^{\bar{r}}_t=A^r_t$ for some $r\in {\Bbb R}'$ such that
$t\in I_r$,

(4) ${\cal S}_{\bar{r}}=\langle S^{\bar{r}}_t:t\in I_{\bar{r}}\rangle$
where $S^{\bar{r}}_t=\bigcup_{s\subseteq t}S_s$ and
$S_s=\bigcup\{S^r_s:(\exists r\in {\Bbb R}')(s\in I_r)\}$.

Notice that from the argument above all $p^{\bar{r}}_t$'s,
$A^{\bar{r}}_t$'s and $S^{\bar{r}}_t$'s are well--defined.
We need to show $\bar{r}\in {\Bbb R}$. It is obvious that $\bar{r}$ is a
common lower bound of all elements in ${\Bbb R}'$ if $\bar{r}\in {\Bbb R}$.

It is easy to see that $\bar{r}$ satisfies 
(1), (2), (3), (4) and (5) in the definition of a condition in
$\Bbb R$. Let's check (6). 

Suppose $t\in I_{\bar{r}}$ and $\alpha\in S^{\bar{r}}_t$.
We want to show that there exists an $a\in (A^{\bar{r}}_t)_{\alpha^{\bar{r}}
_t}$ such that $p^{\bar{r}}_t\forces a\in\dot{f}(\alpha)$.
Let $r\in {\Bbb R}'$ be such that $t\in I_r$, let $r'\in {\Bbb R}'$ and
$s\in I_{r'}$ be such that $s\subseteq t$ and $\alpha\in S^{r'}_s$.
Since $r$ and $r'$ are compatible, then there exists an $r''\in {\Bbb R}$
such that $r''\leqslant r$ and $r''\leqslant r'$. 
By the facts that \[p^{\bar{r}}_t=
p^r_t=p^{r''}_t,\; A^{\bar{r}}_t=A^r_t=A^{r''}_t,\; S^{r'}_s\subseteq 
S^{r''}_s\subseteq S^{r''}_t\] and $r''\in {\Bbb R}$ we have now
that there exists an $a\in (A^{\bar{r}}_t)_{\alpha^{\bar{r}}_t}$ such that
$p^{\bar{r}}_t\forces a\in\dot{f}(\alpha)$.\quad $\square$ (Claim 3.2)

\bigskip

Next we are going to apply $G\!M\!A$ in $M$ to the poset $\Bbb R$ to
construct a $\Bbb P$-name for a Kurepa subtree in $M[G]$.

For each $t\in\omega_1^{<\omega_1}$ define
\[D_t=\{r\in {\Bbb R}:t\in I_r\}.\]
For each $p\in {\Bbb P}$ define
\[E_p=\{r\in {\Bbb R}:(\exists t\in I_r)(p^r_t\leqslant p)\}.\]
For each $\alpha<\omega_1$ define
\[F_{\alpha}=\{r\in {\Bbb R}:(\forall s\in I_r) (\exists t\in I_r) 
(s\subseteq t\mbox{ and }ht(A^r_t)>\alpha)\}.\] For each $\alpha<\kappa$ define
\[O_{\alpha}=\{r\in {\Bbb R}:(\forall s\in I_r)(\exists t\in I_r)
(s\subseteq t\mbox{ and }[\alpha,\kappa)\cap S^r_t\not=\emptyset)\}.\]

\bigskip

{\bf Claim 3.3}\quad All those $D_t$, $E_p$, $F_{\alpha}$ and
$O_{\alpha}$'s are dense in $\Bbb R$.

Proof of Claim 3.3:\quad Let $r_0$ be an arbitrary element in $\Bbb R$.

We show first that for every $t\in\omega_1^{<\omega_1}$ the set
$D_t$ is dense in $\Bbb R$, {\em i.e.} there is an $r\in D_t$ such that 
$r\leqslant r_0$. It's done if $t\in I_{r_0}$. Let's assume that 
$t\not\in I_{r_0}$. Let \[t_0=\bigcup\{s\in I_{r_0}:s\subseteq t\}.\]

Case 1:\quad $t_0\in I_{r_0}$.

Find a sequence $\{p_s:t_0\subseteq s\subseteq t\}$ in $\Bbb P$ such
that $p_{t_0}=p^{r_0}_{t_0}$ and 
\[(\forall s,s')(t_0\subseteq s\subseteq s'\subseteq t\leftrightarrow
p_{s'}\leqslant p_s.\] The sequence $\{p_s:t_0\subseteq s
\subseteq t\}$ exists because $\Bbb P$ is $\omega_1$--closed.
Let \[I_r=I_{r_0}\cup\{s:t_0\subsetneq s\subseteq t\}.\] For any $s\in I_r$,
if $s\in I_{r_0}$, then let \[p^r_s=p^{r_0}_s,\; A^r_s=A^{r_0}_s\mbox{ and }
S^r_s=S^{r_0}_s.\] Otherwise let \[p^r_s=p_s,\; A^r_s=A^{r_0}_{t_0}
\mbox{ and }S^r_s=S^{r_0}_{t_0}.\]
It is easy to see that $r\in D_t$ and $r\leqslant r_0$.

Case 2:\quad $t_0\not\in I_{r_0}$, {\em i.e.} $I_{r_0}$ has no least
element which is above $t$.

Let \[I_r=I_{r_0}\cup\{s:t_0\subseteq s\subseteq t\}.\] Again by
$\omega_1$--closedness we can find \[\{p_s:t_0\subseteq s\subseteq t\}
\subseteq {\Bbb P}\] such that $p_{t_0}$ is a lower bound of \[\{p^{r_0}_s:
s\subseteq t_0\mbox{ and }s\in I_{r_0}\}\] and 
\[(\forall s,s')(t_0\subseteq s\subseteq s'\subseteq t\leftrightarrow
p_{s'}\leqslant p_s).\] Let
\[A'_{t_0}=\bigcup\{A^{r_0}_s:s\in I_{r_0}\mbox{ and }s\subseteq t_0\}\]
and let \[S_{t_0}=\bigcup\{S^{r_0}_s:s\in I_{r_0}
\mbox{ and }s\subseteq t_0\}.\]
If the height of $A'_{t_0}$ is a successor ordinal, then let $A_{t_0}=
A'_{t_0}$. If the height of $A'_{t_0}$ is a limit ordinal, then
we have to add one more level to $A'_{t_0}$.
For any $\beta\in S_{t_0}$ let $s'\subseteq t_0$ and $s'\in I_{r_0}$
be such that $\beta\in S^{r_0}_{s'}$. Then for any $s\in I_{r_0}$
such that $s'\subseteq s\subseteq t_0$ there exists an $a_{s,\beta}\in
(A^{r_0}_s)_{\alpha^{r_0}_s}$ such that $p^{r_0}_s
\forces a_{s,\beta}\in\dot{f}(\beta)$. Now let \[a_{\beta}=\bigcup\{a_{s,\beta}:
s'\subseteq s\subseteq t_0\mbox{ and }s\in I_{r_0}\}\] and let
\[A_{t_0}=A'_{t_0}\cup\{a_{\beta}:\beta\in S_{t_0}\}.\]
It is easy to see that 

(1) the height of $A_{t_0}$ is a successor ordinal,

(2) for every $s\subsetneq t_0$ the tree $A_{t_0}$ is an end--extension
of $A^{r_0}_s$, {\em i.e.}
\[A_{t_0}\!\res\!ht(A^{r_0}_s)=A^{r_0}_s,\]

(3) for every $\beta\in S_{t_0}$ there exists an 
$a_{\beta}$ in the top level of $A_{t_0}$ 
such that \hfill\break $p_{t_0}\forces a_{\beta}\in\dot{f}(\beta)$.

Now for every $s\in I_r$, if $r\in I_{r_0}$, then let \[p^r_s=p^{r_0}_s,\;
A^r_s=A^{r_0}_s\mbox{ and }S^r_s=S^{r_0}_s.\] Otherwise let \[p^r_s=p_s,\;
A^r_s=A_{t_0}\mbox{ and }S^r_s=S_{t_0}.\] It is easy to see that 
$r\in D_t$ and $r\leqslant r_0$.

We show now that for every $p\in {\Bbb P}$ the set $E_p$ is dense in $\Bbb R$.
We want to find an $r\in E_p$ such that $r\leqslant r_0$.
If there exists an $t\in I_{r_0}$ such that $p^{r_0}_t\leqslant p$, then $r_0\in
E_p$. Let's assume that for every $t\in I_{r_0}$ $p^{r_0}_t\not\leqslant p$.
Let \[t_0=\bigcup\{t\in I_{r_0}:p\leqslant p^{r_0}_t\}.\]

Case 1:\quad $t_0\in I_{r_0}$.

Let $t'=t_0\hat{\;}\langle 0\rangle$, {\em i.e.} $t'$ is a successor of $t_0$.
It is clear that $t'\not\in I_{r_0}$.
Let $I_r=I_{r_0}\cup\{t'\}$. For every $t\in I_r$, if $t=t'$, then
let \[p^r_t=p,\; A^r_t=A^{r_0}_{t_0}\mbox{ and }S^r_t=S^{r_0}_{t_0}.\] Otherwise
let \[p^r_t=p^{r_0}_t,\; A^r_t=A^{r_0}_t\mbox{ and }S^r_t=S^{r_0}_t.\]
Then we have $r\in E_p$ and $r\leqslant r_0$.

Case 2:\quad $t_0\not\in I_{r_0}$.

Let $I_r=I_{r_0}\cup\{t_0\}$.
We construct $S_{t_0}$, $A'_{t_0}$ and then $A_{t_0}$ exactly same as
we did in the proof of Case 2 about the denseness of the set $D_t$. For every
$t\in I_r$, if $t=t_0$, then let \[p^r_t=p,\; A^r_t=A_{t_0}\mbox{ and }S^r_t
=S_{t_0}.\] Otherwise let \[p^r_t=p^{r_0}_t,\; A^r_t=A^{r_0}_t\mbox{ and }
S^r_t=S^{r_0}_t.\] Now $r\in E_p$ and $r\leqslant r_0$. Notice also
that $E_p$ is open, {\em i.e.}
\[(\forall p',p''\in {\Bbb P})(p'\leqslant p''\wedge p''\in E_p\rightarrow
p'\in E_p).\]

We show next that for every $\alpha\in\omega_1$ the set $F_{\alpha}$
is dense in $\Bbb R$. We need to find an $r\in F_{\alpha}$ such that
$r\leqslant r_0$.

Let $I_r\supseteq I_{r_0}$ be such that $I_r$ is a countable subtree
of $\omega_1^{<\omega_1}$, $I_r\smallsetminus I_{r_0}$ is an antichain
and for every $s\in I_{r_0}$ there is a $t\in I_r\smallsetminus I_{r_0}$
such that $s\subseteq t$.
For every $t\in I_r\smallsetminus I_{r_0}$ let $p_t\in {\Bbb P}$ be such that
$p_t\leqslant p^{r_0}_s$ for every $s\in I_{r_0}$ and 
$s\subseteq t$, let \[S^r_t=
\bigcup\{S^{r_0}_s:s\in I_{r_0}\mbox{ and }s\subseteq t\}\] and let
\[A'_t=\bigcup\{A^{r_0}_s:s\in I_{r_0}\mbox{ and }s\subseteq t\}.\]
If $ht(A'_t)$ is a successor ordinal, then let $A_t=A'_t$. Otherwise
let \[A_t=A'_t\cup\{a_{\beta}:\beta\in S^r_t\}\] where \[a_{\beta}=
\bigcup\{a\in A'_t: p_t\forces a\in\dot{f}(\beta)\}.\]
Since $S^r_t$ is countable and $\Bbb P$ is $\omega_1$--closed, then
there exists an $p^r_t\leqslant p_t$ such that for every $\beta\in S^r_t$
there exists an $a\in\omega_1^{\alpha}$ such that 
$p^r_t\forces a\in\dot{f}(\beta)$. Let \[A^r_t=A_t\bigcup
\{a\in\omega_1^{\leqslant\alpha}:(\exists\beta\in S^r_t)(p^r_t\forces a\in
\dot{f}(\beta))\}.\] Then $ht(A^r_t)\geqslant\alpha$ is a successor ordinal
and for every $\beta\in S^r_t$ there exists an $a$ in the top level of
$A^r_t$ such that $p^r_t\forces a\in\dot{f}(\beta)$.
For every $t\in I_r\smallsetminus I_{r_0}$ we have already defined $p^r_t$,
$A^r_t$ and $S^r_t$. If $t\in I_{r_0}$, then let \[p^r_t=p^{r_0}_t,\;
A^r_t=A^{r_0}_t\mbox{ and }S^r_t=S^{r_0}_t.\] Hence $r\in F_{\alpha}$
and $r\leqslant r_0$.

We show next that $O_{\alpha}$ for every $\alpha<\kappa$ is dense in $\Bbb R$,
{\em i.e.} finding an $r\in O_{\alpha}$ such that $r\leqslant r_0$.

By imitating the proof of the denseness of $F_{\alpha}$
we can find an $r'\leqslant r_0$ such that $I_{r'}\smallsetminus I_{r_0}$ 
is an antichain and for every $s\in I_{r'}$ there exists an $t\in I_{r'}
\smallsetminus I_{r_0}$ such that $s\subseteq t$. 
For every $t\in I_{r'}\smallsetminus I_{r_0}$ 
fix a $\bar{t}$ which is an successor of $t$
(for example $\bar{t}=t\hat{\;}\langle 0\rangle$).
Let \[I_r=I_{r'}\cup\{\bar{t}:t\in I_{r'}\smallsetminus I_{r_0}\}.\] 
For every $t\in I_{r'}$ let
\[p^r_t=p^{r'}_t,\; A^r_t=A^{r'}_t\mbox{ and }S^r_t=S^{r'}_t.\]
For every $\bar{t}$ with $t\in I_{r'}\smallsetminus I_{r_0}$ 
we want to construct $p^r_{\bar{t}}$,
$A^r_{\bar{t}}$ and $S^r_{\bar{t}}$. If there is a $\beta\in S^{r'}_t$
which is greater than $\alpha$, then let $p^r_{\bar{t}}$ be any proper
extension of $p^r_t$, let $A^r_{\bar{t}}=A^{r'}_t$ and let $S^r_{\bar{t}}=
S^{r'}_t$. Otherwise, first, pick an $a$ in the top level of $A^{r'}_t$,
then choose a $\beta\in\kappa\smallsetminus\alpha$ and a $p\leqslant p^{r'}_t$
such that $p\forces a\in\dot{f}(\beta)$. This can be done because
\[1_{\Bbb P}\forces(\forall t\in\dot{T})(\forall\alpha\in\kappa)(\exists\beta
\in\kappa\smallsetminus\alpha)(t\in\dot{f}(\beta))\] is true in $M$.
Now let \[p^r_{\bar{t}}=p,\; A^r_{\bar{t}}=A^{r'}_t\mbox{ and }S^r_{\bar{t}}
=S^{r'}_t\cup\{\beta\}.\] It is easy to see that $r\in O_{\alpha}$ and
$r\leqslant r_0$.\quad $\square$ (Claim 3.3)

\bigskip

By applying $G\!M\!A$ in $M$ we can find an $\Bbb R$--filter $H$
such that $H\cap D_t\not=\emptyset$, 
$H\cap F_{\alpha}\not=\emptyset$ and $H\cap E_p\cap O_{\alpha'}\not=\emptyset$
for each $t\in\omega_1^{<\omega_1}$, each $\alpha\in\omega_1$, 
each $p\in {\Bbb P}$ and each $\alpha'\in\kappa$.

Since $D_t$ is dense for every $t\in\omega_1^{<\omega_1}$, then
\[I_H=\bigcup\{I_r:r\in H\}=\omega_1^{<\omega_1}.\] Let \[{\Bbb P}_H=
\bigcup\{{\Bbb P}_r:r\in H\}\] and let \[{\cal A}_H=\bigcup\{{\cal A}_r:
r\in H\}.\] Notice that for any $r,r'\in H$ and for any $t\in I_r\cap I_{r'}$
we have $p^r_t=p^{r'}_t$ and $A^r_t=A^{r'}_t$ because $r$ and $r'$
are compatible. So now for every $t\in I_H$ we can define
$p_t=p^r_t$ for some $r\in H$ and define $A_t=A^r_t$ for some $r\in H$.
It is clear that the map $t\mapsto p_t$ is an isomorphism
between $I_H$ and ${\Bbb P}_H$, {\em i.e.} for any $s,t\in I_H$
we have $s\subseteq t$ iff $p_t\leqslant p_s$. It is also clear that
the map $t\mapsto A_t$ is a homomorphism from $I_H$ to ${\cal A}_H$,
{\em i.e.} for any $s,t\in I_H$ we have $s\subseteq t$ implies
$(A_t)\!\res\!ht(A_s)=A_s$.

\bigskip

{\bf Claim 3.4}\quad For each $t\in I_H$ the set 
$\{p_{t\hat{\;}\langle\gamma\rangle}:
\gamma\in\omega_1\}$ is a maximal antichain below $p_t$ in $\Bbb P$.

Proof of Claim 3.4:\quad 
Let $\gamma$ and $\gamma'$ be two ordinals in $\omega_1$. 
Since $I_H=\omega_1^{<\omega_1}$ and $H$ is a filter, there exists an $r\in H$
such that $t\hat{\;}\langle\gamma\rangle$, $t\hat{\;}\langle\gamma'\rangle\in
I_r$. Hence $p^r_{t\hat{\;}\langle\gamma\rangle}$ and $p^r_{t\hat{\;}\langle
\gamma'\rangle}$ are incompatible. So $\{p_{t\hat{\;}\langle\gamma\rangle}:
\gamma\in\omega_1\}$ is an antichain.

Suppose that $p\in {\Bbb P}$ and $p\leqslant p_t$
such that $p$ is incompatible with any of $p_{t\hat{\;}\langle\gamma\rangle}$'s.
Let $r\in H\cap E_p$. Then there is an $s\in I_r$ 
such that $p_s=p^r_s\leqslant p$.
Since $p_s\in {\Bbb P}_H$, then $p_s<p_t$ implies $t\subsetneq s$.
Hence there exists an $\gamma\in\omega_1$ such that $t\hat{\;}\langle\gamma
\rangle\subseteq s$. This means that $p_s\leq 
p_{t\hat{\;}\langle\gamma\rangle}$, {\em i.e.} $p$ and 
$p_{t\hat{\;}\langle\gamma\rangle}$ are compatible, a contradiction.
\quad $\square$ (Claim 3.4)

\bigskip

We now work in $M[G]$. Since $G$ is a $\Bbb P$--generic filter over $M$,
then ${\Bbb P}_H\cap G$ is a linearly ordered subset of ${\Bbb P}_H$.
Let $T_G=\bigcup\{A_t:p_t\in G\}$.

\bigskip

{\bf Claim 3.5}\quad $T_G$ is a Kurepa subtree of $T$ in $M[G]$.

Proof of Claim 3.5:\quad
Since for every $p_t\in G$ we have $p_t\forces A_t\subseteq\dot{T}$, it
is clear that $T_G\subseteq T$ in $M[G]$. For any $p_s,p_t\in G$
we have $p_t\leqslant p_s$ implies $s\subseteq t$ which implies
$(A_t)\!\res\! ht(A_s)=A_s$. Hence $T_G$ is an end--extension
of $A_t$ for every $p_t\in G$. This implies that every level of $T_G$
is a level of some $A_t$, hence is countable.

We want to show now that $T_G$ has at least $\kappa$ branches.
Suppose $|{\cal B}(T_G)|<\kappa$. Then there exists an $\alpha\in\kappa$
such that for every $\beta\in\kappa\smallsetminus\alpha$ the function value
$f(\beta)$ is not a branch of $T_G$. So there is a $p\in {\Bbb P}_H$ 
and there is an $\alpha\in\kappa$ such that \[p\forces (\forall\beta\in\kappa
\smallsetminus\alpha)(\dot{f}(\beta)\mbox{ is not a branch of }
T_{\dot{G}}).\] On the other hand, since $H\cap E_p\cap O_{\alpha}
\not=\emptyset$, then there exists an 
$r\in H\cap O_{\alpha}\cap E_p$.
In $M$ let $s\in I_r$ be such that $p_s\leqslant p$ 
and there is a $\beta\in S^r_s$
such that $\beta>\alpha$. Then for every $t\in I_H$, $s\subseteq t$,
there is an $t'\in I_H$, $t\subseteq t'$, such that 
\[p_{t'}\forces a\in\dot{f}(\beta)\] for some $a\in (A_{t'})_{ht(A_{t'})}$.
This shows that \[p_s\forces\dot{f}(\beta)\mbox{ is a branch of }T_{\dot{G}},\]
which contradicts $p_s\leqslant p$ and \[p\forces (\forall\beta\in\kappa
\smallsetminus\alpha)(\dot{f}(\beta)\mbox{ is not a branch of }T_{\dot{G}}).\]
Hence $T_G$ has at least $\kappa$ branches in $M[G]$.\quad $\square$ 
(Claim 3.5)

\bigskip

Now we conclude that $M[G]\models T\mbox{ has a Kurepa subtree }T_G$.
\quad $\square$

\begin{theorem}

It is consistent with \ch that there exist essential Kurepa trees and
there are no essential Jech--Kunen trees.

\end{theorem}

\noindent {\bf Proof:}\quad
Let $M$ be a model of $C\!H$ and $2^{\omega_1}=\lambda>\omega_2$ plus
$G\!M\!A$. Let ${\Bbb K}_{\lambda}\in M$. Suppose $G_{\lambda}$ is a 
${\Bbb K}_{\lambda}$--generic
filter over $M$. We are going to show that $M[G_{\lambda}]$ is a model of \ch
in which there exist essential Kurepa trees and there are no essential
Jech--Kunen trees.

It is easy to see that $M[G_{\lambda}]$ satisfies \ch. Lemma 1 implies that
there exist essential Kurepa trees. We need only to show that
in $M[G_{\lambda}]$ there are no essential Jech--Kunen trees.

Assume $T$ is a Jech--Kunen tree in $M[G_{\lambda}]$. We need to show that
$T$ has a Kurepa subtree in $M[G_{\lambda}]$. 
Since $|T|=\omega_1$, then there is an
$I\subseteq\lambda$ of cardinality $\omega_1$ in $M$ such that 
$T\in M[G_I]$, where \[G_I=\{p\in G_{\lambda}:dom(l_p)\subseteq I\}.\] 
We claim that \[{\cal B}(T)\cap M[G_{\lambda}]\subseteq M[G_I].\] 
If the claim is true, then $T$ is a Jech--Kunen tree in $M[G_I]$. Suppose that 
$B\in {\cal B}(T)\cap (M[G_{\lambda}]\smallsetminus M[G_I])$. Then there
is a $J\subseteq\lambda\smallsetminus I$ such that $B\in M[G_I][H_J]$
where $H_J$ is a $Fn(J,T_{G_I},\omega_1)$--generic filter over $M[G_I]$.
Let $\dot{B}$ be a $Fn(J,T_{G_I},\omega_1)$--name for $B$. For any 
$J'\subseteq\lambda\smallsetminus (I\cup J)$ such that $|J'|=|J|$
there is an isomorphism $\pi$ from $Fn(J,T_{G_I},\omega_1)$ to
$Fn(J',T_{G_I},\omega_1)$ induced by a bijection between $J$ and $J'$.
Since in $M[G_{\lambda}]$, the branches $(\dot{B})_{H_J}$ and 
$(\pi_*(\dot{B}))_{H_{J'}}$ are different, then $T$ has at least
$\lambda$ branches. This contradicts that $T$ is a Jech--Kunen tree.
Let $T$ have $\delta$ branches
in $M[G_I]$. Since ${\Bbb K}_I$ has size $\omega_1$ and is 
$\omega_1$--closed, then it contains a dense subset which is isomorphic to
${\Bbb P}=(\omega_1^{<\omega_1},\supseteq)$ in $M$. Hence there is a
${\Bbb P}$--generic filter $G$ over $M$ such that $M[G]=M[G_I]$.
By Lemma 3, the tree $T$ has a Kurepa subtree in $M[G]$. Obviously, the Kurepa
subtree is still a Kurepa subtree in $M[G_{\lambda}]$, so $T$ is not
an essential Jech--Kunen tree in $M[G_{\lambda}]$. \quad $\square$

\section{Yes Essential Jech--Kunen Trees, No Essential Kurepa Trees}

In this section we will construct a model of \ch plus the existence of
a thick Kurepa tree, in which there are essential Jech--Kunen trees and 
there are no essential Kurepa trees. The arguments in this section
are a sort of ``symmetric'' to the arguments in the last section.

We first take a model $M$ of $C\!H$ and $2^{\omega_1}=\lambda>\omega_2$
plus a thick Kurepa tree, where $\lambda^{<\lambda}=\lambda$ in $M$, as our
ground model. We then extend $M$ to a model $M[G]$ of $C\!H$ and 
$2^{\omega_1}=\lambda>\omega_2$ plus $G\!M\!A$ by a $\lambda$--stage
iterated forcing (see [B] for the model and forcing). It has been proved
in [Ji1] that in $M[G]$ there are neither essential Jech--Kunen trees nor
essential Kurepa trees. In stead of
taking a model of $G\!M\!A$ as our ground model like we did in \S 1, we
consider this $\lambda$--stage iterated forcing as a part of our
construction because it will be needed later (see also [Ji1, Theorem 5]).
Next we force with a $\omega_1$--closed poset ${\Bbb J}_{S,\kappa}$ in $M[G]$
to create a generic essential Jech--Kunen tree, where $S$ is a 
stationary--costationary subset of $\omega_1$. Again, the hard part is
to prove that forcing with ${\Bbb J}_{S,\kappa}$ over $M[G]$ will not
create any essential Kurepa trees.

\bigskip

Recall that for $T$, a tree, $m(T)$ denote the set
\[\{t\in T:(\forall s\in T)(s\leqslant_T t\rightarrow s=t)\}.\]
Let $I$ be any index set and let $S$ be a subset of $\omega_1$. 
We define a poset ${\Bbb J}_{S,I}$ such that
$p$ is a condition in ${\Bbb J}_{S,I}$ iff $p=(A_p,l_p)$ where

(1) $A_p$ is a countable subtree of $\omega_1^{<\omega_1}$,

(2) $l_p$ is a function from some countable subset of $I$ to $m(A)$.

\noindent For any $p,q\in {\Bbb J}_{S,I}$ define $p\leqslant q$ iff

(1) $A_q\subseteq A_p$,

(2) for every $t\in A_p\smallsetminus A_q$ either there is an 
$s\in m(A_q)$ such that $s\subseteq t$ or
that $\alpha<ht(A_q)$ and $\alpha$ is a limit ordinal imply
\[\alpha=\bigcup\{ht(s):s\in A_q\mbox{ and }s\subseteq t\}\not\in S.\]

(3) $dom(l_q)\subseteq dom(l_p)$ and $(\forall\alpha\in dom(l_q))
(l_q(\alpha)\subseteq l_p(\alpha))$.

\begin{lemma}

($C\!H$) ${\Bbb J}_{S,I}$ is $\omega_1$--closed and  
$\omega_1$--linked.

\end{lemma}

\noindent {\bf Proof:}\quad We show first that ${\Bbb J}_{S,I}$ is 
$\omega_1$--linked. For any $p,q\in {\Bbb J}_{S,I}$, if $A_p=A_q$, then the
condition $(A_p,l_p\cup l_q)$ is a common extension of $p$ and $q$. Because
there are only $\omega_1$ different countable subtrees of 
$\omega_1^{<\omega_1}$, it is clear that ${\Bbb J}_{S,I}$ is the union of
$\omega_1$ linked sets.

We now show that ${\Bbb J}_{S,I}$ is $\omega_1$--closed. 
Let $\{p_n:n\in\omega_1\}$
be a decreasing sequence in ${\Bbb J}_{S,I}$. 
Let $A=\bigcup_{n\in\omega}A_{p_n}$
and let $D=\bigcup_{n\in\omega}dom(l_{p_n})$. For each $i\in D$ let
\[l(i)=\bigcup\{l_{p_n}(i):n\in\omega\mbox{ and }i\in dom(l_{p_n})\}.\]
Define a condition $p\in {\Bbb J}_{S,I}$ such that 
\[A_p=A\cup\{l(i):i\in D\}\mbox{ and }l_p=l.\]
We claim that $p$ is a lower bound of the sequence $\{p_n:n\in
\omega\}$. It suffices to show that for any $n$ and for any 
$t\in A_p\smallsetminus A_{p_n}$ either there exists an 
$s\in m(A_{p_n})$ such that $s\subseteq t$ or
that $\alpha<ht(A_{p_n})$ and $\alpha$ is a limit ordinal imply
\[\alpha=\bigcup\{ht(s):s\in A_{p_n}\mbox{ and }s\subseteq t\}\not\in S.\]
If $t\in A$, then there is an $k>n$ such that $t\in A_{p_k}$.
Hence either there is an $s\in m(A_{p_n})$ such that $s\subseteq t$ 
or that $\alpha<ht(A_{p_n})$ and $\alpha$ is a limit ordinal imply
\[\alpha=\bigcup\{ht(s):s\in A_{p_n}\mbox{ and }s\subseteq t\}\not\in S\]
because $p_k\leqslant p_n$. If $t=l(i)$ for some $i\in D$, then 
because of $t\not\in A_{p_n}$, there is a $k>n$ and 
there is a $t'\in A_{p_k}\smallsetminus A_{p_n}$ such that $t'\subseteq t$.
Hence either there is an $s\in m(A_{p_n})$ such that $s\subseteq t'\subseteq t$
or that $\alpha<ht(A_{p_n})$ and $\alpha$ is a limit ordinal imply
\[\alpha=\bigcup\{ht(s):s\in A_{p_n}\mbox{ and }s\subseteq t'\}\not\in S\]
because $p_k\leqslant p_n$. \quad $\square$

\bigskip

\noindent {\bf Remark:}\quad Again, we may consider the poset ${\Bbb J}_{S,I}$
as a two--step iterated forcing ${\Bbb J}_{S,I'}*Fn(I\smallsetminus I',
T_{\dot{G}_{I'}},\omega_1)$, where $I'$ is a subset of $I$, 
$T_{G_{I'}}=\bigcup\{A_p:p\in G_{I'}\}$ for a generic filter $G_{I'}$ of 
${\Bbb J}_{S,I'}$ and 
$Fn(I\smallsetminus I',T_{\dot{G}_{I'}},\omega_1)$ is a countable
support product of $|I\smallsetminus I'|$--copies of $T_{\dot{G}_{I'}}$.
The map \[p=(A_p,l_p)\mapsto ((A_p,l_p\!\res\! I'),l_p\!\res\! I\smallsetminus 
I')\] is a dense embedding from ${\Bbb J}_{S,I}$ to
${\Bbb J}_{S,I'}*Fn(I\smallsetminus I',T_{\dot{G}_{I'}},\omega_1)$.

\bigskip

We now define $S$--completeness of a tree $T$. Let $\alpha$ be a 
limit ordinal and let $T$ be a tree with $ht(T)=\alpha$. Let
$S$ be a subset of $\alpha$. Then $T$ is called $S$--complete if
for every limit ordinal $\beta\in S$ and every $B\in {\cal B}(T\!\res\!\beta)$
the union $\bigcup B\in T_{\beta}$, {\em i.e.} every strictly decreasing
sequence of $T$ has a greatest lower bound $b$ in $T$ if $ht(b)\in S$.

\begin{lemma}

Let $M$ be a model of $C\!H$ and let 
${\Bbb J}_{S,I}\in M$ where $S\subseteq\omega_1$
and $I$ is an index set in $M$. Suppose $G$ is a ${\Bbb J}_{S,I}$--generic
filter over $M$. Then the tree $T_G=\bigcup_{p\in G}A_p$ is $(\omega_1
\smallsetminus S)$--complete in $M[G]$.

\end{lemma}

\noindent {\bf Proof:}\quad
Let $\alpha\in\omega_1\smallsetminus S$ be a limit ordinal and let
$B$ be a branch of $T_G\!\res\!\alpha$. We need to show that $t=\bigcup B\in
T_G$. The set $B$ is in $M$ because
${\Bbb J}_{S,I}$ is $\omega_1$--closed and $B$ is countable.
Let $p_0\in G$ be such that $B\subseteq A_{p_0}$. It is clear that
\[p_0\forces B\subseteq T_{\dot{G}}.\]
Let \[D_B=\{p\in {\Bbb J}_{S,I}:p\leqslant p_0
\mbox{ and }t=\bigcup B\in A_p\}.\]
Then $D_B$ is dense below $p_0$ because for any $p\leqslant p_0$
the element $p'=(A_p\cup\{\bigcup B\},l_p)$ is a condition
in ${\Bbb J}_{S,I}$ and $p'\leqslant p$ (here we uses the fact that
$\alpha\in\omega_1\smallsetminus S$). Since $p_0\in G$, then
there is a $p\in G\cap D_{B}$. Hence $t=\bigcup B\in T_G$.
\quad $\square$

\begin{lemma}

Let $M$ be a model of $C\!H$. In $M$ let $U$ be a stationary subset of 
$\omega_1$, let $T$ be an $\omega_1$--tree which is $U$--complete and
let $I$ be any index set. Let $K\in M$ be any 
$\omega_1$--tree such that every level of $K$ is countable.
Suppose ${\Bbb P}=Fn(I,T,\omega_1)\in M$ and $G$ is a $\Bbb P$--generic
filter over $M$. Then \[{\cal B}(K)\cap M[G]\subseteq M,\] {\em i.e.}
the forcing adds no new branches of $K$.

\end{lemma}

\noindent {\bf Proof:}\quad
Suppose that $B$ is a branch of $K$ in $M[G]\smallsetminus M$.
Without loss of generality, let's assume that
\[1_{\Bbb P}\forces \dot{B}\in ({\cal B}(K)\smallsetminus M).\]
By a standard argument (see [K2, p. 259]) the statements
\[(\forall p\in {\Bbb P})(\forall\alpha\in\omega_1)
(\exists t\in\omega_1^{\alpha})(\exists p'\leqslant p)(p'\forces t\in\dot{B})\]
and \[(\forall p\in {\Bbb P})(\forall\alpha\in\omega_1)(\forall t\in
\omega_1^{\alpha})(p\forces t\in\dot{B}\longrightarrow\]
\[(\forall\beta\in\omega_1\smallsetminus\alpha)(\exists\gamma\in
\omega_1\smallsetminus\beta)(\exists t_j\in\omega_1^{\gamma})(t_0\not=t_1)
(\exists p_j\leqslant p)(p_j\forces t_j\in\dot{B}))\] 
for $j=0,1$, are true in $M$. 

Let's work in $M$. Let $\theta$ be a large enough cardinal and let
$N$ be a countable elementary submodel of $(H(\theta),\in)$ such that
$K,{\Bbb P},\dot{B}\in N$. Let $\delta=N\cap\omega_1\in U$ (such $N$ exists
because $U$ is stationary). In $M$ we choose an increasing sequence
of ordinals $\{\delta_n:n\in\omega\}$ such that $\bigcup_{n\in\omega}
\delta_n=\delta$. Again in $M$ we construct a set 
\[\{p_s:s\in 2^{<\omega}\}\subseteq {\Bbb P}\cap N\] and a set
\[\{t_s:s\in 2^{<\omega}\}\subseteq K\cap N\] such that

(1) $(\forall s,s'\in 2^{<\omega})(s\subseteq s'\leftrightarrow p_{s'}
\leqslant p_s\leftrightarrow t_{s'}\leqslant t_s)$,

(2) $(\forall s\in 2^{<\omega})(p_s\forces t_s\in\dot{B})$,

(3) $ht(t_s)\geqslant\delta_{|s|}$,

(4) $(\forall i\in dom(p_s))(ht(p_s(i))\geqslant\delta_{|s|})$,

\noindent where $|s|$ means the length of the finite sequence $s$.

Let $p_{\emptyset}=1_{\Bbb P}$ and let $t_{\emptyset}=\emptyset$, the root
of $K$. Assume that we have found $\{p_s:s\in 2^{\leqslant n}\}$ and
$\{t_s:s\in 2^{\leqslant n}\}$ which satisfy (1), (2), (3) and (4) relative
to $2^{\leqslant n}$. Pick any $s\in 2^n$. Since the sentence
\[(\forall p\in {\Bbb P})(\forall\alpha\in\omega_1)(\forall t\in
\omega_1^{\alpha})(p\forces t\in\dot{B}\longrightarrow\]
\[(\forall\beta\in\omega_1\smallsetminus\alpha)(\exists\gamma\in
\omega_1\smallsetminus\beta)(\exists t_j\in\omega_1^{\gamma})(t_0\not=t_1)
(\exists p_j\leqslant p)(p_j\forces t_j\in\dot{B}))\] for $j=0,1$,
is true in $M$, then it is true in $N$. Notice that in $N$ the cardinal
$\omega_1$ is $\delta$. Since $p_s,t_s\in N$, then in $N$ there exist
$p^0,p^1\leqslant p_s$ and there exist $t^0,t^1\in\omega_1^{\gamma}$, 
$t^0\not=t^1$, for some $\gamma\in\delta\smallsetminus\delta_{|s|+1}$
such that \[p^j\forces t^j\in\dot{B}\] for $j=0,1$.
Again in $N$ we can extend $p^0$ and $p^1$ to $p_{s\hat{\;}\langle 0\rangle}$
and $p_{s\hat{\;}\langle 1\rangle}$ respectively so that
\[(\forall i\in dom(p_{s\hat{\;}\langle j\rangle}))(ht(p_{s\hat{\;}\langle
j\rangle}(i))\geqslant\delta_{|s|+1})\] for $j=0,1$. 
Since $T$ is $U$--complete and for every $f\in 2^{\omega}$, for
every $i\in\bigcup_{n\in\omega} dom(p_{f\res n})$ we have
\[\bigcup\{ht(p_{f\res n}(i)):n\in\omega\mbox{ and }i\in dom(p_{f\res n})\}
=\delta\in U,\] then the condition $p_f$ such that $dom(p_f)=
\bigcup_{n\in\omega}dom(p_{f\res n})$ and \[p_f(i)=\bigcup\{p_{f\res n}(i):
n\in\omega\mbox{ and }i\in dom(p_f)\}\] for every $i\in dom(p_f)$
is a lower bound of $\{p_{f\res n}:n\in\omega\}$ in ${\Bbb P}$.
Here we use the fact that $T$ is $U$--complete so that $p_f(i)\in T$
for every $i\in dom(p_f)$.
Let $t_f=\bigcup_{n\in\omega}t_{f\res n}$. Then $ht(t_f)=\delta$.
Since \[p_f\forces t_{f\res n}\in\dot{B}\] for every $n\in\omega$, then
\[p_f\forces t_f\in\dot{B}\cap K_{\delta}.\]
It is easy to see that if $f,f'\in 2^{\omega}$ are different, 
then $t_f$ and $t_{f'}$ are different. Hence $K_{\delta}$ is uncountable,
a contradiction.\quad $\square$

\begin{lemma}

Let $M$ be a model of $C\!H$ and $2^{\omega_1}=\lambda>\omega_2$ and let
${\Bbb J}_{S,\kappa}\in M$ where $\kappa$ is a cardinal in $M$ such that
$\omega_1<\kappa<\lambda$ and $S$ is a stationary subset of $\omega_1$. 
Suppose that $G$ is a ${\Bbb J}_{S,\kappa}$--generic
filter over $M$. Then in $M[G]$ the tree $T_G=\bigcup_{p\in G}A_p$
is an essential Jech--Kunen tree with $\kappa$ branches.

\end{lemma}

\noindent {\bf Proof:}\quad
It is easy to see that $T_G$ is an $\omega_1$ tree. We will divide the 
lemma into two claims.

\bigskip

{\bf Claim 8.1}\quad For every $\xi\in\kappa$ let \[B(\xi)=\bigcup\{l_p(\xi): 
p\in G\mbox{ and }\xi\in dom(l_p)\}.\] Then \[{\cal B}(T_G)=
\{B(\xi):\xi\in\kappa\}\] and for any two different $\xi$ and $\xi'$
in $\kappa$ the branches $B(\xi)$ and $B(\xi')$ are different.

Proof of Claim 8.1:\quad
Since in $M$, for every $\xi\in\kappa$ and for every $\alpha\in\omega_1$
the set \[D_{\xi,\alpha}=\{p\in {\Bbb J}_{S,\kappa}:\xi\in dom(l_p)\mbox{ and }
ht(l_p(\xi))>\alpha\}\] is dense in ${\Bbb J}_{S,\kappa}$, then $B(\xi)$
is a branch of $T_G$. For any two different $\xi,\xi'\in\kappa$
the set \[D_{\xi,\xi'}=\{p\in {\Bbb J}_{S,\kappa}:\xi,\xi'\in dom(l_p)
\mbox{ and }l_p(\xi)\not=l_p(\xi')\}\] is also dense in ${\Bbb J}_{S,\kappa}$.
So the branches $B(\xi)$ and $B(\xi')$ are different.

We now want to show that all branches of $T_G$ in $M[G]$ are exactly
those $B(\xi)$'s. Suppose that in $M[G]$ the tree $T_G$ has a branch
$B$ which is not in the set \[\{B(\xi):\xi\in\kappa\}.\] Without loss of 
generality, let us assume that \[1_{{\Bbb J}_{S,\kappa}}\forces \dot{B}\in
({\cal B}(T_{\dot{G}})\smallsetminus\{\dot{B}(\xi):\xi\in\kappa\}).\] 

Work in $M$. Let $\theta$ be a large enough cardinal and let $N$ be an
elementary submodel of $(H(\theta),\in)$ such that $\kappa,S,\dot{B},
{\cal B}=\{\dot{B}(\xi):\xi\in\kappa\},{\Bbb J}_{S,\kappa}\in N$ and
if $p\in N\cap {\Bbb J}_{S,\kappa}$, then $dom(l_p)\subseteq N$.
Let $\delta=N\cap\omega_1\in S$. In $M$ we choose an increasing sequence
of countable ordinals $\{\delta_n:n\in\omega\}$ such that
$\delta=\bigcup_{n\in\omega}\delta$. We now want to find a decreasing
sequence $\{p_n:n\in\omega\}\subseteq {\Bbb J}_{S,\kappa}\cap N$ such
that $p_0=1_{{\Bbb J}_{S,\kappa}}$ and for each $n\in\omega$

(1) $(\forall\xi\in dom(l_{p_n}))(\exists t\in A_{p_{n+1}})(p_{n+1}
\forces t\in\dot{B}(\xi)\smallsetminus\dot{B})$,

(2) $(\exists\xi\in dom(l_{p_{n+1}}))(\exists t\in A_{p_{n+1}}\smallsetminus
A_{p_n})(ht(t)\geqslant ht(A_{p_n})\mbox{ and }p_{n+1}\forces t\in\dot{B})$,

(3) $ht(A_{p_n})\geqslant\delta_n$.

Assume we have found $\{p_0,p_1,\ldots,p_n\}$. We now work in $N$.
Let \[dom(l_p)=\{\xi_k:k\in\omega\}\] which is an enumeration in $N$. 
Choose $q_0=p_n\geqslant q_1\geqslant
\cdots$ such that for every $k\in\omega_1$ there is a $t\in A_{q_k}$ such
that \[q_k\forces t\in\dot{B}(\xi_k)\smallsetminus\dot{B}.\]
Assume, in $N$, that we have found $\{q_0,q_1,\ldots,q_k\}$. Since the
sentence \[q_k\forces (\exists t\in T_{\dot{G}})(t\in\dot{B}(\xi_k)
\smallsetminus\dot{B})\] is true in $N$ (because it is true in $H(\theta)$
and $\xi_k\in N$), then there is a $t\in(\omega_1^{<\omega_1})^N
=\delta^{<\delta}$ and there is a $q'\leqslant q_k$ such that
\[q'\forces (t\in T_{\dot{G}}\mbox{ and }t\in\dot{B}(\xi_k)\smallsetminus
\dot{B}).\] Since \[q'\forces A_{q'}\subseteq T_{\dot{G}},\]
then there is a $q_{k+1}\leqslant q'$ such that $t\in A_{q_{k+1}}$.
Since $N\models$ ``${\Bbb J}_{S,\kappa}$ is $\omega_1$--closed'' and
$\{q_k:k\in\omega\}$ is constructed in $N$, then there is a $q\in {\Bbb J}_
{S,\kappa}$ in $N$ such that $q$ is a lower bound of $\{q_k:k\in\omega_1\}$.
Let $\alpha=\max\{ht(A_{p_n}),\delta_{n+1}\}$.
Notice that $\alpha\in\delta$ because $p_n\in N$. Since in $N$
\[q\forces\dot{B}\mbox{ is a branch of }T_{\dot{G}},\] then
\[q\forces (\exists t\in (T_{\dot{G}})_{\alpha+1})(t\in\dot{B}).\]
Hence there is a $\bar{q}\leqslant q$ and there is a $t\in 
(\omega_1^{\alpha+1})^N$ such that \[\bar{q}\forces t\in\dot{B}.\]
We can also assume that $t\in A_{\bar{q}}$.

We now go back to $M$ and let $p_{n+1}=\bar{q}$. This finishes the construction
of \hfill\break $\{p_n:n\in\omega\}$.

Let $p\in {\Bbb J}_{S,\kappa}$ be such that \[dom(l_p)=\bigcup_{n\in\omega}
dom(l_{p_n}),\] for every $\xi\in dom(l_p)$\[l_p(\xi)=a_{\xi}=\bigcup
\{l_{p_n}(\xi):n\in\omega\mbox{ and }\xi\in dom(l_{p_n})\}\] and 
\[A_p=(\bigcup_{n\in\omega}A_{p_n})\cup\{a_{\xi}:\xi\in dom(l_p)\}.\]
By the construction of $p_n$'s we have
\[\bigcup\{ht(t):t\in A_p\mbox{ and }p\forces t\in\dot{B}\}=\delta\in S.\]  
Pick any $t\in A_p$. If $t\not=a_{\xi}$ for any $\xi\in dom(l_p)$, then
we can find a $\gamma\in\omega_1$ such that $t\hat{\;}\langle\gamma\rangle
\not\in A_p$. Extend $t\hat{\;}\langle\gamma\rangle$ to 
$\bar{t}\in\omega_1^{delta}$. Define $\bar{p}$ such that 
\[A_{\bar{p}}=A_p\cup\{u:t\subseteq u\subseteq\bar{t}\}\] and $l_{\bar{p}}
=l_p$. If $t=a_{\xi}$ for some $\xi\in dom(l_p)$, then simply extend
$t$ to $b_{\xi}\in\omega_1^{\delta}$ (if $ht(a_{\xi})=\delta$, then
$b_{\xi}=a_{\xi}$). Define $\bar{p}$ such that \[A_{\bar{p}}=A_p\cup
\{u:t\subseteq u\subseteq b_{\xi}\}\] and \[l_{\bar{p}}=(l_p\!\res\! (dom(l_p)
\smallsetminus\{\xi\}))\cup \{(\xi,b_{\xi})\}.\]
It is easy to see that $\bar{p}\leqslant p$ and $ht(A_{\bar{p}})=\delta+1$.
Let \[a=\bigcup\{t\in A_{\bar{p}}:\bar{p}\forces t\in\dot{B}\}.\]
It is also easy to see that for any $q\leqslant\bar{p}$ the element
$a$ is not in $A_q$. Here we use the fact $\delta\in S$, $\delta$ is a limit
ordinal and $ht(A_{\bar{p}})>\delta$. 
Hence \[\bar{p}\forces\dot{B}\cap T_{\dot{G}}
\subseteq\dot{B}\cap A_{\bar{p}}.\] This contradicts that
\[\bar{p}\forces\dot{B}\mbox{ is a branch of }T_{\dot{G}}.\]
\quad $\square$ (Claim 8.1)

\bigskip

{\bf Claim 8.2}\quad $T_{G}$ has no Kurepa subtree in $M[G]$.

Proof of Claim 8.2:\quad Suppose that $T_G$ has a Kurepa subtree $K$
in $M[G]$. Since $|K|=\omega_1$, then there is an $I\subseteq\kappa$
such that $|I|\leqslant\omega_1$ and $K\in M[G_I]$, where
\[G_I=\{p\in G: dom(l_p)\subseteq I\}.\] Notice that $G_I$ is a
${\Bbb J}_{S,I}$--generic filter over $M$. Since ${\Bbb J}_{S,\kappa}$
is forcing equivalent to ${\Bbb J}_{S,I}*Fn(\kappa\smallsetminus I,
T_{G_I},\omega_1)$ and $T_{G_I}$ is $(\omega_1\smallsetminus S)$--complete
in $M[G_I]$ (notice that $S$ is still stationary--costationary), then
by Lemma 7, the set of all branches of $K$ in $M[G_I]$ is same as the set
of all branches of $K$ in $M[G]$. Hence $K$ is a Kurepa tree in $M[G_I]$.
But by Claim 8.1, the tree $T_G=T_{G_I}$ has only $|I|$ branches
in $M[G_I]$ and $K$ is a subtree of $T_G$. Hence $K$ has at most
$\omega_1$ branches in $M[G_I]$. This contradicts that
$K$ is a Kurepa tree in $M[G_I]$.\quad $\square$ 

\begin{lemma}

Let $M$ be a model of $C\!H$ and $2^{\omega_1}=\lambda>\omega_2$ with 
$\lambda^{<\lambda}=\lambda$. In $M$ let $(({\Bbb P}_{\alpha}:\alpha
\leqslant\lambda),(\dot{\Bbb Q}_{\alpha}:\alpha<\lambda))$ be a $\lambda$--stage
iterated forcing notion used in [B] for a model of $G\!M\!A$. Suppose
that $G_{\lambda}$ is a ${\Bbb P}_{\lambda}$--generic filter over $M$.
In $M[G_{\lambda}]$ let ${\Bbb P}=(\omega_1^{<\omega_1},\supseteq)$ and
let $H$ be a $\Bbb P$--generic filter over $M[G_{\lambda}]$. Then
in $M[G_{\lambda}][H]$ there are no essential Kurepa trees.

\end{lemma}

\noindent {\bf Proof:}\quad For any $\alpha<\lambda$ the poset 
${\Bbb P}_{\lambda}$ can be factored to ${\Bbb P}_{\alpha}*{\Bbb P}^{\alpha}$
and $G_{\lambda}$ can also be written as $G_{\alpha}*G^{\alpha}$ such that
$G_{\alpha}$ is a ${\Bbb P}_{\alpha}$--generic filter over $M$ and
$G^{\alpha}$ is a ${\Bbb P}^{\alpha}$--generic filter over $M[G_{\alpha}]$.
Let $T$ be a Kurepa tree in $M[G_{\lambda}][H]$ with $\lambda$ branches. 
Without loss of generality, let's assume 
that for every $t\in T$ there are exactly
$\lambda$ branches of $T$ passing through $t$ in $M[G_{\lambda}][H]$. 
In $M[G_{\lambda}][H]$ let $f:\omega_2\mapsto {\cal B}(T)$ 
be a one to one function such that
for every $t\in T$ and for every $\alpha<\omega_2$ there exists an $\beta\in
\omega_2\smallsetminus\alpha$ such that $t\in f(\beta)$. Notice that
$\omega_2$ here can be replaced by any regular cardinal $\kappa$ satisfying
$\omega_2\leqslant\kappa<\lambda$.
Without loss of generality, let us assume that
\[1_{\Bbb P}\forces (\dot{T}\mbox{ is a Kurepa tree and }
\dot{f}:\omega_2\mapsto {\cal B}(\dot{T})\]\[\mbox{ is a one to one function
such that }(\forall t\in\dot{T})(\forall\alpha\in\omega_2)(\exists\beta\in
\omega_2\smallsetminus\alpha)(t\in\dot{f}(\beta))).\]
We want now to construct a poset ${\Bbb R}'$ in $M$ such that a filter
$\bar{G}$ of ${\Bbb R}'$ obtained by applying a forcing argument similar to
$G\!M\!A$ in $M[G_{\lambda}]$ will
give us an $\Bbb P$--name for a Jech--Kunen subtree of 
$T$ in $M[G_{\lambda}][H]$.

Let $r$ be a condition in ${\Bbb R}'$ iff $r=(I_r,{\Bbb P}_r,{\cal A}_r,
{\cal S}_r)$ where $I_r$ is a countable subtree of $(\omega_1^{<\omega_1},
\supseteq)$, ${\Bbb P}_r=\langle p^r_t:t\in I_r\rangle$, ${\cal A}_r=
\langle A^r_t:t\in I_r\rangle$ and ${\cal S}_r=\langle S^r_t:t\in I_r\rangle$
such that

(1) ${\Bbb P}_r\subseteq {\Bbb P}$, and for every $t\in I_r$ the element
$A^r_t$ is a nonempty countable subtree of 
$(\omega_1^{<\omega_1},\supseteq)$ of
height $\alpha^r_t+1$ (we will use some $A^r_t$'s to generate a 
Jech--Kunen subtree of $T$) and $S^r_t$ is a nonempty
countable subset of $\omega_2$, (the requirement ``$S^r_t\subseteq\omega_2$''
makes ${\Bbb R}'$ different from ${\Bbb R}$ defined in Lemma 3,)

(2) $(\forall s,t\in I_r)(s\subseteq t\leftrightarrow p^r_t\leqslant p^r_s)$,

(3) $(\forall s,t\in I_r)(s\subseteq t\rightarrow 
A^r_t\!\res\! ht(A^r_s)=A^r_s)$, 

(4) $(\forall s,t\in I_r)(s\subseteq t\rightarrow S^r_s\subseteq S^r_t)$,

(5) $(\forall t\in I_r)(p^r_t\forces A^r_t\subseteq\dot{T})$,

(6) $(\forall t\in I_r)(\forall\alpha\in S^r_t)(\exists a\in 
(A^r_t)_{\alpha^r_t})(p^r_t\forces a\in\dot{f}(\alpha))$.

For any $r,r'\in {\Bbb R}'$, let
$r\leqslant r'$ iff $I_{r'}\subseteq I_r$, and for every $t\in I_{r'}$
\[p^{r'}_t=p^r_t,\; A^{r'}_t=A^r_t\mbox{ and }S^{r'}_t\subseteq S^r_t.\]

\bigskip

{\bf Claim 9.1}\quad The poset ${\Bbb R}'$ is $\omega_1$--linked.

Proof of Claim 9.1:\quad Same as the proof of Claim 3.1.\quad $\square$ 
(Claim 9.1)

\bigskip

{\bf Claim 9.2}\quad The poset ${\Bbb R}'$ is countably compact.

Proof of Claim 9.2:\quad Same as the proof of Claim 3.2.
\quad $\square$ (Claim 9.2)

\bigskip

For each $t\in\omega_1^{<\omega_1}$ define
\[D_t=\{r\in {\Bbb R}':t\in I_r\}.\]
For each $p\in {\Bbb P}$ define
\[E_p=\{r\in {\Bbb R}':(\exists t\in I_r)(p^r_t\leqslant p)\}.\]
For each $\alpha<\omega_1$ define
\[F_{\alpha}=\{r\in {\Bbb R}':(\forall s\in I_r) (\exists t\in I_r) 
(ht(A^r_t)>\alpha)\}.\] For each $\alpha<\omega_2$ define
\[O_{\alpha}=\{r\in {\Bbb R}':(\forall s\in I_r)(\exists t\in I_r)
(s\subseteq t\mbox{ and }[\alpha,\omega_2)\cap S^r_t\not=\emptyset)\}.\]

\bigskip

{\bf Claim 9.3}\quad All those $D_t$, $E_p$, $F_{\alpha}$ and
$O_{\alpha}$'s are dense in ${\Bbb R}'$.

Proof of Claim 9.3:\quad Same as the proof of Claim 3.3.
\quad $\square$ (Claim 9.3)

\bigskip

Note that $|{\Bbb R}'|=\omega_2$. Also note that $M[G_{\lambda}][H]=
M[H][G_{\lambda}]$ because ${\Bbb P}_{\lambda}$ is $\omega_1$--closed.
By the construction of ${\Bbb P}_{\lambda}$ there exists an 
$\beta<\lambda$ such that those dense sets 
$D_t$, $E_p$, $F_{\alpha}$ and $O_{\alpha}$ are 
in $M[G_{\beta}]$, the tree $T$ is in $M[G_{\beta}][H]$ or $\dot{T}$ is in
$M[G_{\beta}]$ and
\[1_{{\Bbb P}_{\beta}}\forces\dot{\Bbb Q}_{\beta}={\Bbb R}',\] {\em i.e.} 
${\Bbb R}'$ is the poset used in $\beta$--th step forcing in
the $\lambda$--stage iteration.

Let $H_{\beta}$ be a ${\Bbb Q}_{\beta}$--generic filter over $M[G_{\beta}]$
such that $G_{\beta}*H_{\beta}=G_{\beta+1}$.

Since $D_t$ is dense for every $t\in\omega_1^{<\omega_1}$, then
\[I_{H_{\beta}}=\bigcup\{I_r:r\in H\}=\omega_1^{<\omega_1}.\] 
Let \[{\Bbb P}_{H_{\beta}}=
\bigcup\{{\Bbb P}_r:r\in H_{\beta}\}\] and let 
\[{\cal A}_{H_{\beta}}=\bigcup\{{\cal A}_r:
r\in H_{\beta}\}.\] 
Notice that for any $r,r'\in H_{\beta}$ and for any $t\in I_r\cap I_{r'}$
we have $p^r_t=p^{r'}_t$ and $A^r_t=A^{r'}_t$ because $r$ and $r'$
are compatible. So now for every $t\in I_{H_{\beta}}$ we can define
$p_t=p^r_t$ for some $r\in H_{\beta}$ and define $A_t=A^r_t$ 
for some $r\in H_{\beta}$.
It is clear that the map $t\mapsto p_t$ is an isomorphism
between $I_{H_{\beta}}$ and ${\Bbb P}_{H_{\beta}}$, 
{\em i.e.} for any $s,t\in I_{H_{\beta}}$
we have $s\subseteq t$ iff $p_t\leqslant p_s$. It is also clear that
the map $t\mapsto A_t$ is a homomorphism from $I_{H_{\beta}}$ to 
${\cal A}_{H_{\beta}}$,
{\em i.e.} for any $s,t\in I_{H_{\beta}}$ we have $s\subseteq t$ implies
$(A_t)\!\res\!ht(A_s)=A_s$.

\bigskip

{\bf Claim 9.4}\quad For each $t\in I_{H_{\beta}}$ the set 
$\{p_{t\hat{\;}\langle\gamma\rangle}:
\gamma\in\omega_1\}$ is a maximal antichain below $p_t$ in $\Bbb P$.

Proof of Claim 9.4:\quad Same as the proof of Claim 3.4.

\bigskip

The next claim is something different from Lemma 3.
Let $T_H=\{A_t:p_t\in H\}$ where $H$ is the ${\Bbb P}$--generic filter
over $M[G_{\lambda}]$.

\bigskip

{\bf Claim 9.5}\quad $T_H$ is a Jech--Kunen subtree of $T$ in 
$M[G_{\lambda}][H]$.

Proof of Claim 9.5:\quad
By the proof of Claim 3.5, we know that $T_H$ is a subtree of $T$
with more than $\omega_1$ branches. It suffices to show that
$T_H$ has exactly $\omega_2$ branches.

Suppose that $T_H$ has more than $\omega_2$ branches then there is a
branch $B$ in $M[G_{\lambda}][H]$ which is not in the range of the function
$f$. Without loss of generality, let's assume that \[1_{\Bbb P}\forces
(\forall \alpha\in\omega_2)(\dot{B}\not=\dot{f}(\alpha)).\]
Let $\dot{B}$ be a $\Bbb P$--name for $B$ and let
\[D_{\dot{B}}=\{r\in {\Bbb R}': (\forall s\in I_r)(\exists t\in I_r)
(s\subseteq t\mbox{ and } ht(\dot{B}\cap A^r_t)<ht(A^r_t))\}.\]
Since $M[G_{\lambda}][H]=M[G_{\beta}][H][G^{\beta}]$ and ${\Bbb P}^{\beta}$
is $\omega_1$--closed in $M[G_{\beta}][H]$, then $B$ is in $M[G_{\beta}][H]$
because any $\omega_1$--closed forcing will not add any new branches to
the Kurepa tree $T$. We assume also that the $\Bbb P$--name $\dot{B}$
is in $M[G_{\beta}]$. Hence the set $D_{\dot{B}}$ is in $M[G_{\beta}]$.
Let \[E_{\dot{B}}=\{p^r_t\in {\Bbb P}_{H_{\beta}}:r\in D_{\dot{B}}\cap H_{\beta}
\mbox{ and }p^r_t\forces ht(\dot{B}\cap A^r_t)<ht(A^r_t)\}.\]

\bigskip

{\bf Subclaim 9.5.1}\quad $D_{\dot{B}}$ is dense in ${\Bbb R}'$.

Proof of Claim 9.5.1:\quad
Let $r_0$ be any element in ${\Bbb R}'$. It suffices to show that there
is an element $r$ in $D_{\dot{B}}$ such that $r\leqslant r_0$. 
Let's first extend $r_0$ to $r'$ such that for every $s\in I_{r_0}$ there
is a $t\in m(I_{r'})$ such that $s\subseteq t$.
Let $t\in m(I_{r'})$. For every $\alpha\in S^{r'}_t$ let $a_{\alpha}
\in (A^{r'}_t)_{\alpha^{r'}_t}$ such that $p^{r'}_t\forces a_{\alpha}
\in\dot{f}(\alpha)$. Since we have
\[p^{r'}_t\forces(\exists u\in\dot{T})(u\in\dot{f}(\alpha)\smallsetminus
\dot{B})\] and $\Bbb P$ is $\omega_1$--closed, then there is
a $u_{\alpha}\supseteq a_{\alpha}$ in $\omega_1^{<\omega_1}$ 
for every $\alpha\in S^{r'}_t$ 
and a $p_t\leqslant p^{r'}_t$ such that for every $\alpha\in S^{r'}_t$
\[p_t\forces u_{\alpha}\in\dot{f}(\alpha)\smallsetminus\dot{B}.\]
Without loss of generality, we can assume that there is
a $\gamma\in\omega_1$ such that $ht(u_{\alpha})=\gamma$
and \[p_t\forces\dot{B}\mbox{ differs from 
all }\dot{f}(\alpha)\mbox{ below }\gamma\]
for every $\alpha\in S^{r'}_t$. 
Let \[I_r=I_{r'}\cup\{\bar{t}:\bar{t}\mbox{ is a successor of }t\mbox{ for }
t\in m(I_{r'})\}.\] For every $t\in I_{r'}$ let \[p^r_t
=p^{r'}_t,\; A^r_t=A^{r'}_t\mbox{ and }S^r_t=S^{r'}_t.\] 
For every $\bar{t}\in I_r\smallsetminus I_{r'}$ let \[p^r_{\bar{t}}=p_t,\; 
A^r_{\bar{t}}=A^{r'}_t\cup\{s:s\subseteq u_{\alpha}
\mbox{ for some }\alpha\in S^{r'}_t\}\mbox{ and }S^r_{\bar{t}}=S^{r'}_t.\]
Now it is easy to see that $r\leqslant r_0$ and $r\in D_{\dot{B}}$.
\quad $\square$ (Subclaim 9.5.1)

\bigskip

{\bf Subclaim 9.5.2}\quad $E_{\dot{B}}$ is dense in ${\Bbb P}_{H_{\beta}}$.

Proof of Subclaim 9.5.2:\quad
Let $p_0\in {\Bbb P}_{H_{\beta}}$. We need to show that there is a
$p\in {\Bbb P}_{H_{\beta}}$ such that $p\leqslant p_0$ and $p\in E_{\dot{B}}$.

Since $p_0\in {\Bbb P}_{H_{\beta}}$, then there is an $r\in H_{\beta}$
such that $p_0=p^r_s$. Since $D_{\dot{B}}$ is dense and $r\in H_{\beta}$,
then there is an $r'\leqslant r$ such that $r'\in H_{\beta}\cap D_{\dot{B}}$.
Since $p^r_s=p^{r'}_s$ and $r'\in D_{\dot{B}}$, then there is a $t\in
I_{r'}$ such that $s\subseteq t$ and
\[p^{r'}_t\forces ht(\dot{B}\cap A^{r'}_t)<ht(A^{r'}_t).\]
Hence we have $p^{r'}_t\leqslant p^r_s=p_0$ and $p^{r'}_t\in E_{\dot{B}}$.
\quad $\square$ (Claim 9.5)

\bigskip

Now the lemma follows because if $B$ of is a branch of $T$, which is not
in the range of $f$, then it is not a branch of $T_{H_{\beta}}$ because
there is an $\alpha\in\omega_1$ and there is an $A=T_{H_{\beta}}\!\res\!\alpha$
such that $ht(B\cap A)<ht(A)$.\quad $\square$

\begin{theorem}

It is consistent with \ch plus the existence of a thick Kurepa tree that
there exist essential Jech--Kunen trees and there are no essential 
Kurepa trees.

\end{theorem}

\noindent {\bf Proof:}\quad
Let $M$ be a model of $C\!H$ and $2^{\omega_1}=\lambda>\omega_2$ such that
in $M$, $\lambda^{<\lambda}=\lambda$ and there is a thick Kurepa tree. Such
model exists by Lemma 1. In $M$ let 
\[(({\Bbb P}_{\alpha}:\alpha\leqslant\lambda),
(\dot{\Bbb Q}_{\alpha}:\alpha<\lambda))\] be the $\lambda$--stage iterated
forcing notion used in [B] for a model of $G\!M\!A$. Suppose $G_{\lambda}$
is a ${\Bbb P}_{\lambda}$--generic filter over $M$. Then \[M[G_{\lambda}]
\models C\!H+2^{\omega_1}=\lambda>\omega_2+G\!M\!A.\]
In $M[G_{\lambda}]$ let $\kappa$ be a cardinal such that $\omega_2\leqslant
\kappa<\lambda$ and let $S$ be a stationary--costationary
subset of $\omega_1$. Suppose that $H$ is a ${\Bbb J}_{S,\kappa}$--generic 
filter over $M[G_{\lambda}]$. Then by Lemma 8, the tree $T_H=
\bigcup\{A_p:p\in H\}$ is an essential Jech--Kunen tree in $M[G_{\lambda}][H]$.
It is obvious that the thick Kurepa trees in $M$ are still thick Kurepa trees
in $M[G_{\lambda}][H]$. We need only to show that there are no essential
Kurepa trees in $M[G_{\lambda}][H]$.

Suppose that $K$ is an essential Kurepa tree in $M[G_{\lambda}][H]$.
Since $|K|=\omega_1$, then there exists an $I\subseteq\kappa$ such that
$|I|=\omega_1$ and $K\in M[G_{\lambda}][H_I]$, where \[H_I=H\cap {\Bbb J}_{S,I}
=\{p\in H: dom(l_p)\subseteq I\}.\] Since ${\Bbb J}_{S,\kappa}$
is forcing equivalent to \[{\Bbb J}_{S,I}*Fn(\kappa\smallsetminus I,
T_{\dot{H}_I},\omega_1))\] and by Lemma 6, the tree $T_{H_I}$ is
$(\omega_1\smallsetminus S)$--complete,
then by Lemma 7, there are no new branches of $K$
in $M[G_{\lambda}][H]$ which are not in $M[G_{\lambda}][H_I]$.
So $K$ is still a Kurepa tree in $M[G_{\lambda}][H_I]$. But the poset
${\Bbb J}_{S,I}$ is $\omega_1$--closed and has size $\omega_1$.
So by Lemma 2, the poset ${\Bbb J}_{S,I}$ is forcing equivalent to
$(\omega_1^{<\omega_1},\supseteq)$. Hence by Lemma 9, 
the Kurepa tree $K$ has a Jech--Kunen subtree
$K'$ in $M[G_{\lambda}][H_I]$. Since every branch of $K'$
is a branch of $K$ and the set of branches of $K$
keeps same in $M[G_{\lambda}][H_I]$ and in $M[G_{\lambda}][H]$, then
$K'$ is still a Jech--Kunen subtree of $K$ in $M[G_{\lambda}][H]$. This
contradicts that $K$ is an essential Kurepa tree in $M[G_{\lambda}][H]$
\quad $\square$

\bigskip

\noindent {\bf Remark:}\quad It is quite easy to build a model of
\ch in which there exist essential Jech--Kunen trees and there are no 
essential Kurepa tree without requiring the existence of a thick
Kurepa tree. Let $M$ be a model of $G\!C\!H$. First, increase $2^{\omega_1}$
to $\omega_3$ by an $\omega_1$--closed Cohen forcing. Then, force with
the poset ${\Bbb J}_{S,\omega_2}$. In the resulting model
$C\!H$ and $2^{\omega_1}=\omega_3$ hold and there is an essential
Jech--Kunen tree. It can be shown easily that there are no
thick Kurepa trees in the resulting model. Hence it is trivially
true that there are no essential Kurepa trees in that model.

\section{New Proofs of Two Old Results.}

In [SJ1], we proved that, assuming the consistency of an
inaccessible cardinal, it is consistent with \ch that there exist
Jech--Kunen trees and there are no Kurepa trees. The model for that 
is constructed by taking Kunen's model for non--existence of Jech--Kunen
trees as our ground model and then forcing with a countable support
product of $\omega_2$ copies of a ``carefully pruned'' tree $T$. 
The way that the tree $T$ is pruned guarantees that 
(1) the forcing is $\omega$--distributive,
(2) forcing does not add any Kurepa trees, (3) $T$ becomes a Jech--Kunen tree
in the resulting model.
In [Ji3], this pruning technique was also used to construct a model of
\ch in which there exist essential Kurepa trees and there exist essential
Jech--Kunen trees. Here we realize that the Jech--Kunen tree obtained by
forcing with that carefully pruned tree in [SJ1] and [Ji3] can be replaced by
a generic Jech--Kunen tree obtained by forcing with ${\Bbb J}_{S,\kappa}$,
the poset defined in \S 2. So now we can reprove those two results
in [SJ1] and [Ji3] without going through a long and tedious construction
of a ``carefully pruned'' tree.

Let $Lv(\kappa,\omega_1)$, the countable support L\'{e}vy collapsing order, 
denote a poset defined by letting $p\in Lv(\kappa,\omega_1)$ iff
$p$ is a function from some countable subset of $\kappa\times\omega_1$
to $2$ such that $p(\xi,\eta)\in\xi$ for every $(\xi,\eta)\in dom(p)$ and
orderd by reverse inclusion.

Let $Fn(\lambda,2,\omega_1)$, the countable support Cohen forcing,
denote a poset defined by letting $p\in Fn(\lambda,2,\omega_1)$ iff
$p$ is a function from some countable subset of $\lambda$ to $2$ and
ordered by reverse inclusion.

\begin{theorem}

Let $\kappa$ and $\lambda$ be two cardinals in a model $M$ such that
$\kappa$ is strongly inaccessible and $\lambda>\kappa$ is regular in $M$.
Let $S\in M$ be a stationary--costationary subset of $\omega_1$ and
let ${\Bbb J}_{S,\kappa}\in M$ be the poset defined in \S 2.  Let
$Lv(\kappa,\omega_1)$ and $Fn(\lambda,2,\omega_1)$ be in $M$. Suppose that
$G\times H\times F$ is a $(Lv(\kappa,\omega_1)\times Fn(\lambda,2,
\omega_1)\times {\Bbb J}_{S,\kappa})$--generic filter over $M$. Then
$M[G][H][F]\models (C\!H+2^{\omega_1}>\omega_2\;+$ there exist
Jech--Kunen trees $+$ there are no Kurepa trees $)$.

\end{theorem}

\noindent {\bf Proof:}\quad
It is easy to see that \[M[G][H][F]\models (C\!H+2^{\omega_1}=\lambda>
\kappa=\omega_2).\] It is also easy to see that $\omega_1$ 
and all cardinals greater than or equal to $\kappa$ in $M$ are preserved.
By Lemma 8, the tree $T_F=\bigcup_{p\in F}A_p$ is a Jech--Kunen tree.
We now need only to show that there are no Kurepa trees in $M[G][H][F]$.
Suppose that $K$ is a Kurepa tree in $M[G][H][F]$. Since $|K|=\omega_1$,
then there exists an $I\subseteq\kappa$ with $|I|=\omega_1$ such that
$K\in M[G][H][F_I]$ where $F_I=F\cap {\Bbb J}_{S,I}$ (recall that
the poset ${\Bbb J}_{S,\kappa}$ is forcing equivalent to ${\Bbb J}_{S,I}*
Fn(\kappa\smallsetminus I,T_{\dot{F}_I},\omega_1)$).
By Lemma 7, the tree $K$ is still a Kurepa tree in $M[G][H][F_I]$.
Since the poset ${\Bbb J}_{S,I}$ is $\omega_1$--closed and has size
$\omega_1$, then by Lemma 2, ${\Bbb J}_{S,I}$ is forcing equivalent
to $Fn(\omega_1,2,\omega_1)$. By a standard argument we
know that $Fn(\lambda,2,\omega_1)\times Fn(\omega_1,2,\omega_1)$ is isomorphic
to $Fn(\lambda,2,\omega_1)$. Hence there is a $Fn(\lambda,2,\omega_1)$--generic
filter $H'$ over $M[G]$ such that $M[G][H][F_I]=M[G][H']$.
But it is easy to see that in $M[G][H']$ there are neither Kurepa trees nor
Jech--Kunen trees. So we have a contradiction that $K$ is a Kurepa tree
in $M[G][H']$.\quad $\square$

\begin{theorem}

Let $M$ be a model of $G\!C\!H$. Let $\kappa$ and $\lambda$ be two
regular cardinals in $M$ such that $\lambda>\kappa>\omega_1$ and let
$S$ be a stationary subset of $\omega_1$ in $M$.
In $M$ let ${\Bbb K}_{\lambda}$ and ${\Bbb J}_{S,\kappa}$ be two posets
defined in \S 1 and \S 2, respectively. Suppose that $G\times H$ is
a ${\Bbb K}_{\lambda}\times {\Bbb J}_{S,\kappa}$--generic filter over
$M$. Then \[M[G\times H]\models (C\!H+2^{\omega_1}=\lambda>\kappa>\omega_1\;+\]
\[\mbox{ there exist essential Kurepa trees }+\mbox{ there exist 
essential Jech--Kunen trees }).\]

\end{theorem}

\noindent {\bf Proof:}\quad
It is easy to see that $M[G\times H]$ is a model of $C\!H$ and $2^{\omega_1}=
\lambda>\kappa>\omega_1$.
Since ${\Bbb K}_{\lambda}$ and ${\Bbb J}_{S,\kappa}$ are $\omega_1$--closed,
then ${\Bbb K}_{\lambda}$ is absolute with respect to $M$, and $M[H]$ and
${\Bbb J}_{S,\kappa}$ is absolute with respect to $M$ and $M[G]$.
By Lemma 8, the tree $T_H=\bigcup_{p\in H}A_p$ is an essential Jech--Kunen 
tree in $M[G][H]$. By Lemma 1, the tree $T_G=\bigcup_{p\in G}A_p$
is an essential Kurepa tree because $M[G][H]=M[H][G]$.
\quad $\square$

\bigskip

Department of Mathematics

University of California

Berkeley, CA 94720

{\em e-mail: jin@@math.berkeley.edu}

\bigskip

Institute of Mathematics, 

The Hebrew University, 

Jerusalem, Israel.

\bigskip

Department of Mathematics, 

Rutgers University, 

New Brunswick, NJ, 08903, USA.

\bigskip

{\em Sorting:} The first address is the first author's and the last
two are the second author's. 

\end{document}